\documentstyle{amsppt}
\magnification=1200

\def\ss{\smallskip}
\def\ms{\medskip}

\def\pa{\partial}

\def\oonL{\mathop{{\hbox{$L^{m,2}$}\kern -20pt\raise7pt
\hbox{$\circ$}}}}

\def\O{\Omega}
\def\o{\omega}

\def\bigtimes{\mathop{{\lower2.95pt\hbox{$\wedge$}\kern-6.666pt\raise2.95pt
\hbox{$\vee$}}}\limits}
\def\crd{\cr\noalign{\vskip4pt}}

%??

\def\eq{\eqalign}

\def\e{\epsilon}

\def\R{${\text{\bf R}}^N$}
\def\no{\noindent}
\def\bs{\bigskip}

\def\g{\gamma}

\def\b{\beta}
\def\a{\alpha}

\def\da{d_{\partial \O}(x)}
\def\d{\da}

\let\hacek=\v% "Sparar" check-accent som /hacek, innan /v omdefinieras.
\def\v{\vert}

%\magnification 1200
\hsize=5in
\vsize=7.6in
\baselineskip=16pt

\def\Wmps0O{W^{m,p}_0(\O,\d^s)}
\def\Wmpt0O{W^{m,p}_0(\O,\d^t)}

\NoRunningHeads
\widestnumber\key{ABCDEFG}
\NoBlackBoxes
\nologo

\topmatter
\title
Removability:\ \ \ \ \ \ \ \ \ \ \ \ \ \ \ \ \ \ \ \ \ \ \ \ \ \ \ \ 
\ \ \ \ \ \ \ \ \ \ \ \ 
\\
First order Sobolev Space, \ \ \ \ \ \ \ \ \ \ \ \ \ \ \ \ \ \ \ \ \ \ 
\\
PDE solutions, holomorphic functions.\ \ \ \ \ \ \ 
\\
-- Sobolev space constructed in a new way.
\endtitle
\author
Andreas Wannebo
\endauthor
\address
Department of Mathematics
\newline The Royal Institute of Technology
\newline Sweden
\endaddress
\abstract
Conditions for removability for sets with respect to first order Sobolev space
functions, holomorphic functions and other solutions to first order
PDE:s -- also nonlinear -- are proved.
A new kind of representation formula for first order Sobolev space is proved. 
This result can also be used as a new definition of these spaces and
is useful for the discussion of removability of sets in \R 
with respect to to this space.

The problem of removability for first order Sobolev space gives
answers to a question on the classification of Julia sets, see [DOA-HUB].

This treatment of Sobolev space removability is a follow up of a paper by
Peter Jones, see [JON]. 
\endabstract
\endtopmatter
\document

\heading 1. Introduction and Discussion 
\endheading

First it should be said that this paper is a follow up of a paper by
Peter Jones, [JON], and that the main contribution here is Theorem 2.16,
which gives a new way to define first order Sobolev space, more about
this later in this section.

The main part of the present paper have been given as work document,
[WAN], and this have been sent to Peter Jones some years ago. 
This have had effect on his later work on the problem.
Peter Jones and Stanislav Smirnov have recently produced a paper, [JON-SMI].
\bs

After this has been said we begin a longer introduction/discussion 
to give background and motivation.
\bs

In a general way of speaking a removability problem can be formulated as
follows. Let $\Cal F$ and $\Cal G$ be classes of functions from say \R\ to say 
$\text{\bf R}$ or $\text{\bf C}$.
Suppose that this class $\Cal G$ is a nice class e.g. the 
holomorphic functions on the Riemann sphere $\text{\bf C}$ or the 
real-valued Sobolev functions in $W^{m,p}(\text{\bf R}^N)$.

Suppose that $\Cal F$ has the same kind of description 
as that of $\Cal G$ but with a deficit in domain so to speak, 
e.g. are holomorphic on the complement of a set
$S$ or are Sobolev functions on the complement of $S$.
Then suppose there are extra conditions as well -- this is often
the case.

The question of removability in these settings is to find
conditions that implies that the set $S$ can be removed, i.e.
all functions in $\Cal F$ are restrictions of functions in 
$\Cal G$. Negative answers are of course interesting as well.

Next we go through some removability problems.
For the moment let's assume Sobolev space to be known with some
definition or the other. Proper definitions will be given later.
\bs

{\bf Notation 1.1.} Let $\O$ be an open subset of \R\ and
let $W^{1,p}(\O)$ denote Sobolev space of order one with exponent $p$
and with domain $\O$.
\bs

In order to get nicer formulas we also introduce an adjusted notation.
\bs

{\bf Notation 1.2.} 
Let $\O$ be open in \R.
Denote with $W^{1,p}(\O,\text{\bf R}^N)$ the set of functions 
from \R\ to $\text{\bf R}$, which have their restriction to $\O$ 
in $W^{1,p}(\O)$.
\bs

Now we choose to discuss removability in terms of compacts rather then 
closed sets in order to conform with the formulations in [JON].

Let $C(\text{\bf R}^N)$ denote the space of real continuous functions on
\R.
\bs

{\bf Definition 1.3.} Let $p\ge 1$, let $K$ be compact in $\text{\bf R}^N$.
Then put $K=\O^c$. That $K$ has the $J_p$-property if the
following two subspaces are equal

$$
C(\text{\bf R}^N)\cap W^{1,p}(\O,\text{\bf R}^N)
=C(\text{\bf R}^N)\cap W^{1,p}(\text{\bf R}^N).
\leqno(1.1)
$$

Note that the inclusion RHS$\subset$LHS in (1.1) is obvious and
(1.1) can be reduced to LHS$\subset$RHS.
\bs

Hence that a compact has the $J_p$-property has the meaning that the compact is
removable in $W^{1,p}$ with continuity as extra property.
\bs

Next we introduce some removability problems involving
holomorphic functions.

For some of these there is a background involving capacity concepts. 
There are two standard capacities for analytic functions,
analytic capacity $\g(K)$ and continuous analytic capacity $\a(K)$. 
For definitions and results concerning these see the book by Garnett, [GAR].
\bs

Theorem. If $\o$ is an open subset of $\text{\bf C}$ and 
compact $K$ subset of $\o$ then the property that
all $u$ continuous in the $\o$ and holomorphic in $\o\setminus K$
also are holomorphic in the whole $\o$ holds if and only if $\a(K)=0$.
\bs

We can give some of the contents of this theorem without giving the
definition of $\a(K)$.
\bs

Theorem. Besicovitch, see [BES]. If $K$ is compact in $\text{\bf C}$
and if $K$ is a a countable union of sets, each with finite 
one-dimensional Hausdorff measure, then $\a(K)=0$. 
\bs

Hence $\a(K)=0$ implies that $K$ is removable in the sense given.
\bs

Next we give another removability problem for holomorphic functions.
\bs

{\bf Definition 1.4.} 
Let $K$ be a compact subset of the Riemann sphere $\text{\bf C}$. The
set $K$ is said to have the $HR$-property (Holomorphically Removable)
if that $F$ is a homeomorphism from $\text{\bf C}$ to $\text{\bf C}$ with 
$F$ holomorphic off $K$ implies that $F$ is a M\"obius transformation.
\bs

That this is a removability property follows from that if $F$ is a
M\"obius transformation then $K$ was removable and on the other hand
if $K$ is removable then $F$ is defined everywhere as an analytic
function on the Riemann sphere and then that $F$ also is a homeomorphism 
implies that $F$ is a M\"obius transformation.
\bs

Hence if $K$ is removable in the sense related to continuous analytic capacity
then $K$ has the $HR$-property. Hence also $\a(K)=0$ implies that $K$ has the 
$HR$-property.
\bs

One background for the interest in the HR-property is found in 
the classical work by Doady and Hubbard, [DOA-HUB], 
on iteration of one complex variable quadratic polynomials.

Doady and Hubbard proved among many things that the corresponding Julia sets
can be divided into disjoint classes in a natural way and that
such a class of Julia sets in this classification has 
exactly one member if and only if any Julia set in the class has the 
HR-property.
\bs

We follow Jones and identify $\text{\bf C}$ and $\text{\bf R}^2$, see
[JON]. This way the compacts can be seen as subset of both and
the HR-property and the $J_p$-property can be compared.
\bs

The results by Peter Jones in [JON] can be summed up as follows.
Here (i) is an observation and the main part is (ii).
\ss

\roster
\item"(i)" 
The $J_2$-property implies the HR-property. 
\endroster

\roster
\item "(ii)"
The boundary of a John domain has the $J_2$-property.
\endroster
\ss

The definition of John domains is well-known. We refer to [JON].

This result may seem strange at first, since it is obvious that very
simple domains, e.g. a domain that form a simple cusp, have the
$J_2$-property and they are clearly not John domains.

The explanation is that Peter Jones was posed exactly the question 
to prove the HR-property for the boundary of a general John domain
in order to give an answer to a question concerning a special situation in
Dynamical System theory.
\bs

However from the function space point of view the correct setting of
this removability problem 
is with the weaker property quasicontinuity instead of continuity.
The quasicontinuity property will be defined in the next section.

We define other kinds of removability akin to the previous ones.

It is natural to think about these different definitions of 
removability  and what implications there are between them.
\newpage

\roster
\item "(i)"
The continuous $J_p$-property. This is written $J_p(cont)$.
\endroster

\roster
\item "(ii)"
The quasicontinuous $J_p$-property. Written $J_p(quasicont)$ or
just simply $J_p$.
\endroster

\roster
\item "(iii)"
The $J_p$ and $J_p(cont)$ properties locally in the sense
of $L^p$-norms is written $J_{p,loc}$ and $J_{p,loc}(cont)$.
\endroster

\roster
\item "(iv)"
It follows that $J_p$ implies $J_p(cont)$.
\endroster

\roster
\item "(v)"
It follows that $J_{p,loc}=J_p$ and $J_{p,loc}(cont)=J_p(cont)$
by the use of partition of unity.
\endroster

\roster
\item "(vi)"
$HL_p$, $HL_p(cont)$, $HL_{p,loc}$, etc. 
defined as removability of holomorphic functions with respect to the 
$L^p$-norm of the gradient and otherwise in accordance with the
definition of $J_p$ etc.
\endroster
\ms

{\bf Some open questions.}  
\ms

\roster
\item "{\bf I:}"
Is there a compact $K$ such that $J_p$ does not hold
for $K$ but $J_p(cont)$ holds.
\endroster

\roster
\item "{\bf II:}"
Is there a compact $K$ such that
$J_p$, ($J_p(cont)$), does not hold for $K$, but holds for $J_q$,
($J_q(cont)$), $p\not=q$.
\endroster

\roster
\item "{\bf III:}"
Jones: Is there a compact $K$ such that
$HR$ holds for $K$ but not for $J_2(cont)$.
\endroster
\bs

{\bf A discussion on results in general.}
\ss

A problem is not solved just because there exists a proved necessary
and sufficient condition. The real issue at stake is how many of the
interesting cases are covered and can be tested.
There is here also a question about what type of descriptions, geometrical,
topological etc. that can be covered.
\ss

-- These or similar kinds of questions are also essential for other
fields of problems.

It would be profitable if such questions could be more (publicly) discussed.
\bs

In our case there are plenty of results on removability though it is
for the moment not possible to reach as far into the tough cases as 
wanted, compare here [JON] and [JON-SMI].
\bs

Since $J_p(cont)$ is simpler to study the main effort is spent on the
more difficult $J_p$ case. 
\bs

-- It is interesting that by posing a question in this function space-wise 
correct way we are lead to make a new, alternative and interesting definition 
of first order Sobolev space.
\bs

Next we formulate a well-known fact, which will be important. Since
the proof is short it is included.
\bs

{\bf Theorem 1.5.} Let $u$ be a function on \R, let $u$ be 
absolutely continuous
on almost all lines in the $x_1$-direction. Let $u, D_1u \in L^1_{loc}$,
then the weak derivative $D_1u$ exists and is given as the pointwise 
derivative $D_1u$.
\bs

{\bf Proof.} On each of these special lines the derivative $D_1u$
exists a.e. by the absolute continuity, but this is a.a. lines,
hence the pointwise derivative $D_1u$ is defined a.e. in the domain.

However the pointwise derivative can be identified with the
weak derivative a distribution.

Simply take a test function $\varphi$ and evaluate the action of
$D_1u$ on $\varphi$ as usual.

$$
-\int uD_1\varphi dx=-\int\int_LuD_1\varphi dtdx'
=\int\int_LD_1u\varphi dtdx',
\leqno(1.2)
$$
by partial integration, which works since $u$ is absolute continuous 
on allmost all lines $L$.

End of Proof.
\bs

{\bf Definition 1.6.} Let $V^{1,p}(\O)$ be what could be called the
rough Sobolev space on the domain $\O$ in \R. This definition is that
both the function and gradient is pointwise defined a.e. Hence the 
corresponding Sobolev space norm makes sense. It should be finite as well.
\bs

-- In the same way as is done in Notation 1.2 we assume that the functions in
$V^{1,p}(\text{\R})$ are defined in all \R.
\bs

{\bf The main result, Theorem 2.16.}
\ss

There are two equivalent formulations. Here we discuss (2.11).

In words (2.11) states that the rough Sobolev space $V^{1,p}(\O)$ 
when intersected with the the set of quasicontinuous functions gives
a direct sum in the vector space sense. 

The first summand is in fact the ordinary Sobolev space but with
functions defined in a precise sense quasi everywhere.
The second summand contains the the singular part of this subspace
of the rough Sobolev space also in a precise sense.
\bs

How does this apply to the $J_p$-property?
\bs

Let the open set $\O$ in the formula be $K^c$. Let $u\in W^{1,p}(K^c)$.
Now it is not difficult to see that a necessary condition for
$J_p$ to hold is that the Lebesgue measure of $K$ is zero. 
-- Hence assume this. 
It follows then that $u\in V^{1,p}(\text{\bf R}^N)$. 
\bs

Then study the $J_p$-property or the $J_p(cont)$-property.

-- In either case the property is true if and only if 
the second direct summand is the empty set.
\bs

{\bf Corollary 1.7.} The $J_p$-property, ($J_p(cont)$-\-property),
for comp\-act $K$ in \R\, holds 
if and only if there does not exist a function 
$f\in Q\cap V^{1,p}(\text{\bf R}^N)$,
($f\in \text{\bf C}\cap V^{1,p}(\text{\bf R}^N)$),
such that $f\in W^{1,p}(K^c)$ and $f$ is singular on a set of
parallell lines of nonzero $(N-1)$-dim Lebesgue measure.
\bs

{\bf Corollary 1.8.} It is neccessary for the $J_p$ property to be untrue that
there is a set of parallell lines of nonzero $(N-1)$-dim Lebesgue measure that
intersect $K$ uncountably many times, since a singular function on a line
always has uncountable support.
\bs

{\bf Corollary 1.9.} It is sufficient for the $J_p$ property to hold that
there are $N$ non-degenerate directions such that the lines in these directions
intersect $K$ with a $(N-1)$-dim Lebesgue measures for these lines are equal
to zero.
\bs

{\bf Removability with respect to PDO:s.}
\ss

Given a function $u$ which is the weak solution to a first order PDE
in $\text{\R}\setminus S$. Where $S$ is closed and of Lebesgue measure zero. 
Suppose that $D_1u$ and $u$ are in $L^1$. Then by Theorem 1.5 if
a.a. lines in this direction intersect $S$ countably many times then
$D_1u$ is a global weak derivative. This can be done in any direction.
Also say linear complex combinations can be proven to give say
$\bar\pa u$ as global weak derivative.

Observe. The PDO does not have be linear above.
\bs

{\bf Remark.} Hence $J_2(cont)$ implies $HL_2(cont)$, which 
implies $HR$. (Cf. area theorem.)

Though we knew this before from [JON].
\bs
\bs

\heading
2. Main section
\endheading

To begin with some notations and definitions.
\bs

{\bf Definition 2.1.} Let $\O$ be an open set in \R and let $T(\O)$ 
be the set of measurable functions from $\O$ to {\bf R}, then the  
$L^p$-norm of a function in $u\in T(\O)$ is

$$
||u||_{L^p(\O)}=(\int_{\O}|u|^pdx)^{1/p}.
\leqno(2.1)
$$
\bs

The corresponding first order Sobolev space norm is defined as follows.
\bs

{\bf Definition 2.2.}
$$
||u||_{W^{1,p}(\O)}=\sum_{|\a|\le 1}||D^{\a}u||_{L^p(\O)},
\leqno(2.2)
$$
where $\a$ is a multiindex. It has always to be clear what definition of 
derivative that is used, since this can change.
\bs

{\bf Comment 2.3.} There are other definitions of first order Sobolev space
norms. They are usually equivalent (in the usual analysis sense)
to the one given above.
\bs

Our way to get at the promised new point on first order
Sobolev space is by several definitions of similar spaces, different 
super- and subspaces. Hence many more definitions will be needed.
\bs

First we give the Sobolev definition of Sobolev space.
It was however not the first one. It is an early exampel of the use 
of weak derivatives long before elaborate distribution theory.
\bs

{\bf Definition 2.4.} Let $\O$ be open subset of \R, then
$$
W^{1,p}(\O)=
\{u\in L^p(\O): ||u||_{W^{1,p}(\O)}<\infty\},
\leqno(2.3)
$$
where derivatives are taken in the weak sense.
\bs

{\bf Definition 2.5.} Let $\O$ be open subset of \R, then
the rough Sobolev space is defined as

$$
V^{1,p}(\O)=
\{u\in L^p(\O): ||u||_{W^{1,p}(\O)}<\infty\}.
\leqno(2.4)
$$

The derivatives are taken in the usual sense and  exist a.e.
\bs

We assume the concept of absolute continuity to be familiar.
\bs

{\bf Definition 2.6.} Denote by $[AC]_j$ the subset of $T(\O)$ of members
that they are absolutely continuous with respect to almost all
$(N-1)$-dimensional Lebesgue measure) lines parallell with the $x_j$
axis in any closed cube parallel to the axii and subset of $\O$.
Furthermore denote

$$
[AC]=
\bigcap_{j=1}^N[AC]_j
\leqno(2.5)
$$

The following are important and in fact show how
much Sobolev space is related to $L^p$-space and this not only formally.
\bs

{\bf Definition 2.7.} Let $u$ be a function from $\O$ open in \R\ to {\bf R}.
Let $\Cal F$ be a family of subsets of $\O$ thta contains the open subsets.
Let $\Psi$ a mapping from $\Cal F$ to the nonnegative real numbers.
The family ${\Cal F}$ is now denoted ${\Cal F}_{\Psi}$.

The function $u$ is said to be quasi continuous on $\O$ with respect
to $\Psi$, if to every $\e>0$ there exists an open $\o$ subset of
$\O$ such that $\Psi(\o)<\e$ and $u|_{\O\setminus \o}$ is continuous.
\bs

It is well-known that in Sobolev space the place of Lebesgue measure
is taken by a set function called capacity i.e. in this case a very special
capacity. Often the definition of this capacity is taken from
potential theory and more specifically non-linear potential theory,
i.e. when $p\not=2$. 
\bs

We will use the intrinsic definition of capacity that comes
from the definition of the Sobolev space. 
For $1<p$ these definitions give equivalent values.
\bs

The capacity $C_{1,p}$, also written as $(1,p)$-capacity, is a set 
function mapping a certain family of sets into the nonnegative reals.

C.f. $\Psi$ above.

It is closely connected to the structure of the first order Sobolev space.
\bs

{\bf Definition 2.8.} Let $K$ be a compact set in \R. Then define
$$
C_{1,p}(K)=\inf\ \{\sum_{|\a|\le 1}||D^{\a}\varphi||_
{L^p(\text{\bf R}^N)}:
\varphi\in C^{\infty}_0(\text{\bf R}^N), 
\ \varphi|_K\ge 1\}.
\leqno(2.6)
$$
\bs

The following steps serves the purpose to generalize the domain of
definition for the capacity.
\bs

{\bf Definition 2.9.} Let $\o$ be an open set in \R. Then define
$$
C_{1,p}(\o)=\sup_K\ \{C_{1,p}(K): K\ {\text{compact}},\ K\subset \o\}
\leqno(2.7)
$$
\bs

Now we do a presentation proper, though there is a possibility to
make a short-cut in the presentation here.
\bs

{\bf Definition 2.10.} A set $M\in\text{\bf R}^N$ is capacitable if

$$
\inf\{C_{1,p}(\o):\o\ \text{open},\ M\subset \o\}=
\sup\{C_{1,p}(K):K\ \text{compact},\ K\subset M\}\}
\leqno(2.8)
$$

The family of all capacitable subsets of \R\ is denoted ${\Cal F}_C$.
\bs

Hence we let capacity be a mapping from ${\Cal F}_C$ to the
nonnegative reals and with the above procedure as definition.
\bs

The capacity defined here has a many properties. We will recall what
needed, when neeeded.

Yet we say the following. In the first order case it follows by so
called truncation that it is enough to consider nonnegative
$\varphi$:s in the first step 
in the definition of $(1,p)$-capacity above. It is then obvious
that this capacity is subadditive, cf. measures that are additive.

For higher order Sobolev space this question is dealt with non-linear
potential theory.
\bs

Next we have some facts, well-known from non-linear potential theory.
This theory assumes in our case that $p>1$.
\bs

{\bf Fact 2.11.}
The space $W^{1,p}(\text{\bf R}^N)$ equals a disjunct union of sets of 
functions, which are equal almost everywhere within each such set. 
There is in each set at least one member, which is quasi continuous 
with respect to $C_{1,p}$.
\bs

{\bf Definition 2.12.}
$$
Q=\ (Q^{1,p}(\O))\ =\{u\in T(\O):
u\text{ is quasi continous w.r.t. }C_{1,p}\}
\leqno (2.9)
$$

Furthermore:
\bs

{\bf Fact 2.13.}
Uniqueness: Let $u,u'\in W^{1,p}(\text{\bf R}^N)$ and say that
$u=u'$ almost every where and say furthermore that $u$ and $u'$ are
quasi continuous (with respect $(1,p)$-capacity) then $u=u'$ quasi
everywhere, i.e. they differ on a set of $(1,p)$-capacity zero only.
\bs

{\bf Definition 2.14.} Let $S_j$ be the subset of elements $T(\O)$ with $j$th partial
derivative zero, as defined a.e.
\bs

{\bf Definition 2.15.} Let $\text{\hacek X}_j$ be the subset $T(\O)$ 
such that the elements are independent of $x_j$ on almost all lines
parallell to the $j$th axis.
\bs

The new point of view of (first order) Sobolev space 
is expressed by the equalities (2.10) and (2.11) below. It serves also as
an alternative definition of this space.
\bs

{\bf Theorem 2.16.} Let $1\le p$. 
Let the functions below be considered equal if they are
equal quasi everywhere, i.e. differ only on a set of (1,p)-capacity zero. 
The quotients and direct sum are here taken in the vector space sense,
i.e. no topology concerns. Then it holds that

$$
Q\cap V^{1,p}(\O)/\sum_{j=1}^N\ [S_j\cap Q\cap V^{1,p}(\O)/
\text{\hacek X}_j\cap Q\cap V^{1,p}(\O)]=
[AC]\cap Q\cap V^{1,p}(\O)
\leqno(2.10)
$$
or equivalently
$$
\ \ \ Q\cap V^{1,p}(\O)=[AC]\cap Q\cap V^{1,p}(\O)\oplus
\sum_{j=1}^N\ [S_j\cap Q\cap V^{1,p}(\O)/
\text{\hacek X}_j\cap Q\cap V^{1,p}(\O)].
\leqno(2.11)
$$
\bs

{\bf Proof.} It is given as a sequence of lemmas etc.
\bs

{\bf Lemma 2.17.} Observe that from the vector space structure we have
by Zorn's lemma a vector space basis $B$ of $Q\cap V^{1,p}(\O)$.
Since 
$$
 \text{\hacek X}_j\cap Q\cap V^{1,p}(\O))\subset S_j\cap Q\cap V^{1,p}(\O)
\subset Q\cap V^{1,p}(\O)
\leqno(2.12)
$$
as subspaces, there furthermore exists such a basis $B$ that
there are $B_j\subset B$, with the following
properties.
\bs

(i) Every $B_j$ is a basis of
$$
S_j\cap Q\cap V^{1,p}(\O)/\text{\hacek X}_j\cap Q\cap V^{1,p}(\O).
\leqno(2.13)
$$
\bs

(ii) Let $J$ be a subset of the whole numbers $1,2 ... N$, then 
$\bigcup_{j\in J}B_j$ is a basis of 
$$
\sum_{j\in J}\ [S_j\cap Q\cap V^{1,p}(\O)/
\text{\hacek X}_j\cap Q\cap V^{1,p}(\O).
\leqno(2.14)
$$

(iii) With $J$ as above, $\bigcap_{j\in J}B_j$ is a basis of 
$$
\bigcap_{j\in J}S_j\cap Q\cap V^{1,p}(\O)/
\text{\hacek X}_j\cap Q\cap V^{1,p}(\O).
\leqno(2.15)
$$
\bs

The proof the statements (i)--(iii) follows simply by the following
steps. First a
partion of the original vector space in a Venn diagram way and an
ordering of the corresponding ``atomic'' subspaces. Next
by finitely many succesive extensions of basis we get the properties needed.
\bs

It now follows from Lemma 2.17 that if such a basis is chosen then
(2.10) and (2.11) are equivalent just by set operations on the elements
in the basis.
\bs

Next we give some facts that will be needed in order to proceed.
\bs

{\bf Fact 2.18.} The capacity $C_{1,p}(M)\ge m(M)$, where $m$ is Lebesgue
measure.
\bs

Proof: This follows from the definitions. Lebesgue measure $m$ can be
defined as the capacity $C_{0,p}$.
\bs

The possibility of the above definition of Lebesgue measure in fact 
is responsible for the aforementioned similarities between $L^p$ and
Sobolev space.
\bs

{\bf Fact 2.19.} Let $C_{1,p}$ be given in \R. Then the corresponding capacity
value after a projection been made to hyperplane (or any plane) is 
less than or to equal the first capacity value up to a constant 
factor.
\bs

-- This is a result from non-linear potential theory.
\bs

We will assume the definition of Hausdorff measure to be known, at least
as in our case when the defining function $h_{\a}$ is the $\a$-power of 
the cubes side length. 

We here give one of the well-known connections between non-linear 
capacity and Hausdorff measure, proved by non-linear potential theory.
\bs

{\bf Fact 2.20.} (A special case.) Let $M\subset\text{\bf R}^n$ and 
$C_{1,p}(M)=0$, then for $\e>0$ it holds that $h_{n-p+\e}=0$.
\bs

The fact below is immediate from the definitions.

{\bf Fact 2.21.} Let $M$ be as above, then
\ss

(i) $\a<\b$ and $h_{\a}(M)=0$ then $h_{\b}(M)=0$;
\ss

(ii) $h_n$ and $n$-dimensional Lebesgue measure are the same.
\bs

We assume the concept singular to be known. The fact below is well-know
and elementary.
\bs

{\bf Fact 2.22.} A continuous function on a line is the unique sum of an 
absolutely continuous part and a singular part, i.e. direct sum in
vector sense.
\bs

We proceed with the proof and choose to prove (2.10).
\bs

Since there always is the possibility of a  partition of unity 
using smooth functions, it is clear that we can confine ourselves to the
question for an arbitrary cube $B$, such that its closure lies in
$\O$. 
\bs

Hence let $u\in Q\cap V^{1,p}(B)$. Furthermore let $L_j$ be the set of
lines in the $x_j$-direction with the property that $u$ is continuous
on each line in the set $L_j$. To state it clearly we use a
restriction notation for classes of functions with respect to $L_j$
and similar sets, but this shall be interpreted on an individual
function basis, or formally as $L_j$ etc. as dependent on the function
involved. 

-- This is a minor change versus the normal notations.
\bs

{\bf Notation 2.23.}
Let
$$
\eq{
E_j&=\bigcup_{l\in L_j}l
\crd
E&=\cap_jE_j
\cr}
\leqno(2.16)
$$
\bs

Now we want to show that 
$$
m_{N-1}(E_j^c\cap M_{N-1})=0
\leqno(2.17)
$$

Here $M_{N-1}$ is a $N-1$ dimensional hyperplan perpendicular to the
$x_j$-direction and $m_{N-1}$ is the $N-1$ dimensional Lebesgue measure.

To show this we start with that $u\in Q$. Then the points where $u$
is discontinuous can be kept in an open set $\o_{\e}$ such that
$C_{1,p}(\o_{\e})\le \e$. The set $D\o_{\e}$ can be projected
orthogonally to $M_{N-1}$ and we get the set 
$\o_{\e}^{M_{N-1}}$.
\bs

Now by Fact 2.19 we have 
$$
C_{1,p}(\o_{\e}^{M_{N-1}})\le A\e
\leqno(2.18)
$$. 

However we if go to the limit we really have 
$$
C_{1,p}(E_j^c\cap M_{N-1})=0.
\leqno(2.19)
$$

Now by Fact 2.20 we have that for $\e'>0$ it holds that
$$
h_{N-p+\e'}(E_j^c\cap M_{N-1})=0.
\leqno(2.20)
$$

This is rewritten as
$$
h_{N-1-(p-1)+\e'}(E_j^c\cap M_{N-1})=0.
\leqno(2.21)
$$

Now if $\e'$ is chosen less than $p-1$ -- or why not equal 
to -- which is possible since we  as always assume $p>1$, then we have
$$
h_{N-1}(E_j^c\cap M_{N-1})=0.
\leqno(2.22)
$$

But by Fact 2.21 this implies that
$$
m_{N-1}(E_j^c\cap M_{N-1})=0,
$$
i.e. (2.17).
\bs

Now we are in a position to discuss (2.10) with domain the cube $B$.
The LHS of (20) can be written as follows, as a consequence of Lemma
2.17,

$$
\eq{
Q\cap &V^{1,p}(B)/\sum_{j=1}^N\ [S_j\cap Q\cap V^{1,p}(B)/
\text{\hacek X}_j\cap Q\cap V^{1,p}(B)]
\crd
&=
\bigcap_{j=1}^N\ (Q\cap V^{1,p}(B)/[S_j\cap Q\cap V^{1,p}(B)/
\text{\hacek X}_j\cap Q\cap V^{1,p}(B)]).
\cr}
\leqno(2.23)
$$

However we should rather use some set $L_j'$ instead of $L_j$, since
the ``almost every'' prefix in the definitions of $\text{\hacek X}_j$ etc.
Take the intersection(s) and get $L_j'$. 

This set obviously
contains almost all lines in $x_j$-direction as well -- a finite union
of sets of measure zero has measure zero. Hence this adjustment has no
significance.

Furthermore redefine $E_j$ to $E_j'$ and $E$ to $E'$.
\bs

Our next step is to set out to prove what is needed but with restriction
to $E'$, i.e.
\bs

{\bf Step 2.24.} It holds 
$$
\bigcap_{j=1}^N\ (Q\cap V^{1,p}(B)/[S_j\cap Q\cap V^{1,p}(B)/
\text{\hacek X}_j\cap Q\cap V^{1,p}(B)]|_{E'})=
$$
$$=[AC]\cap Q\cap V^{1,p}(B)|_{E'}.
\leqno(2.24)
$$

Now again we will return to the basis argument via Zorn's lemma and
Venn diagramms used to prove Lemma 2.17.  In this case we get the
following
$$
\bigcap_{j=1}^N Q\cap V^{1,p}(B)/[S_j\cap Q\cap V^{1,p}(B)/
\text{\hacek X}_j\cap Q\cap V^{1,p}(B)]|_{E'}=
$$
$$
=V^{1,p}(B)\cap\bigcap_{j=1}^N\  Q/[S_j\cap Q/
\text{\hacek X}_j\cap Q]|_{E'}.
\leqno(2.25)
$$

However it follows from the definitions of $S_j$ and $\text{\hacek X}_j$ as
well as the unique decomposition of a contionuous function on a line
as a sum of an absolutely continuous part and a singular part 
(Fact 2.22), that
$$
Q/[S_j\cap Q/\text{\hacek X}_j\cap Q]|_{E'_j}=[AC]_j\cap Q|_{E_j'}.
\leqno(2.26)
$$

Now the information adds up. We use (2.24) and (2.25) together with the
definition of $[AC]$ and get Step 2.24, i.e. (2.24).
\bs

This result Step 2.24 has because of the very definition of $[AC]$ 
the implication that a drop of the restriction gives
$$
\eq{
\bigcap_{j=1}^N&\ (Q\cap V^{1,p}(B)/[S_j\cap Q\cap V^{1,p}(B)/
\text{\hacek X}_j\cap Q\cap V^{1,p}(B)]\subset
\crd
&\subset 
[AC]\cap Q\cap V^{1,p}(B).
\cr}
\leqno(2.27)
$$
\bs

{\bf Fact 2.25.} It holds that

$$
[AC]\cap V^{1,p}(\O)\subset W^{1,p}(\O);\text{ Here }\O=B
\leqno(2.28)
$$
\bs

The LHS above was the original definition of first order Sobolev
space. After Sobolev made his definition it has been more in fashion.

-- There are other definitions still.
\bs

Now there is lining up of the inclusions

$$
\eq{
Q\cap& W^{1,p}(B)\subset 
Q\cap V^{1,p}(B)/ \sum_{j=1}^N
[S_j\cap Q\cap V^{1,p}(B)/ \text{\hacek X}_j\cap Q\cap V^{1,p}(B)]
\crd
&\subset
Q\cap [AC]\cap V^{1,p}(B)
\subset 
Q\cap W^{1,p}(B).
\cr}
\leqno(2.29)
$$

Since the begin equals the end, all inclusions actually are equalities.
\bs

Hence we have proved (2.10), i.e. the theorem.
\bs

End of proof.

\enddocument